%an AMS Tex file 
%eightteen pages of output
%of a paper by 
%  
%Dominique Foata and Doron Zeilberger
%
%"A combinatorial proof of Bass's Evaluations of 
%the Ihara-Selberg Zeta Function for Graphs"

\input amstex
%-----------------------------------------------------------------------------
% Beginning of tran.sty
%-----------------------------------------------------------------------------

%  First input the preprint style
%%% ====================================================================
%%% @AMSTeX-style-file{
%%%   filename  = "ams-j.sty",
%%%   version   = "2.1b",
%%%   date      = "1995/04/03",
%%%   time      = "10:06:15 EDT",
%%%   checksum  = "57047 483 915 16049",
%%%   author    = "American Mathematical Society",
%%%   address   = "PO Box 6248, Providence, RI 02940, USA",
%%%   email     = "tech-support@math.ams.org (Internet)",
%%%   supported = "yes",
%%%   keywords  = "",
%%%   abstract  = "This is an AMSTeX documentstyle. It uses the amsppt
%%%     documentstyle as a base and adds a few refinements to give
%%%     printed documents the visual form used for final publication of
%%%     AMS journal articles. There are also documentstyles for
%%%     individual AMS journals that input this documentstyle and add
%%%     journal-specific information (ISSN, journal name, etc).",
%%%   docstring = "The checksum field above contains: CRC-16
%%%     checksum, number of lines, number of words, and number of
%%%     characters, as produced by Robert Solovay's checksum utility.",
%%% }
%%% ====================================================================

%  First input the preprint style
\def\next{AMSPPT}\ifx\styname\next \else\input amsppt.sty \relax\fi

\catcode`\@=11

\indenti=3pc
\pagewidth{30pc} \pageheight{46.5pc}

\font@\fourteenbf=cmbx10 scaled \magstep2

\newtoks\sevenpoint@
\def\sevenpoint{\normalbaselineskip9\p@
 \textonlyfont@\rm\sevenrm \textonlyfont@\it\sevenit
 \textonlyfont@\sl\sevensl \textonlyfont@\bf\sevenbf
 \textonlyfont@\smc\sevensmc \textonlyfont@\tt\seventt
  \textfont\z@\sevenrm \scriptfont\z@\sixrm
       \scriptscriptfont\z@\fiverm
  \textfont\@ne\seveni \scriptfont\@ne\sixi
       \scriptscriptfont\@ne\fivei
  \textfont\tw@\sevensy \scriptfont\tw@\sixsy
       \scriptscriptfont\tw@\fivesy
  \textfont\thr@@\sevenex \scriptfont\thr@@\sevenex
   \scriptscriptfont\thr@@\sevenex
  \textfont\itfam\sevenit \scriptfont\itfam\sevenit
   \scriptscriptfont\itfam\sevenit
  \textfont\bffam\sevenbf \scriptfont\bffam\sixbf
   \scriptscriptfont\bffam\fivebf
 \setbox\strutbox\hbox{\vrule height7\p@ depth3\p@ width\z@}%
 \setbox\strutbox@\hbox{\raise.5\normallineskiplimit\vbox{%
   \kern-\normallineskiplimit\copy\strutbox}}%
 \setbox\z@\vbox{\hbox{$($}\kern\z@}\bigsize@1.2\ht\z@
 \normalbaselines\sevenrm\dotsspace@1.5mu\ex@.2326ex\jot3\ex@
 \the\sevenpoint@}

%       define  a logo for the upper left-hand corner

\def\jourlogo{%
\def\rightheadline{\vtop to 0pt{%
        \sixrm \baselineskip7pt
        \parindent0pt \frenchspacing
AMS JOURNAL STYLE\hfil\break
        Volume \cvol@, \cvolyear@\par\vss}}}

%       macros to be put into the \topmatter for the logo

\def\cvol#1{\gdef\cvol@{\ignorespaces#1\unskip}}
\def\cvolno#1{\gdef\cvolno@{\ignorespaces#1\unskip}}
\def\cmonth#1{\gdef\cmonth@{\ignorespaces#1\unskip}}
\def\cvolyear#1{\gdef\cvolyear@{\ignorespaces#1\unskip}}
\def\cyear#1{\gdef\cyear@{\ignorespaces#1\unskip}\cyear@@#100000\end@}
\def\cpgs#1{\gdef\cpgs@{\ignorespaces#1\unskip}}

\def\cyear@@#1#2#3#4#5\end@{\gdef\cyearmodc@{#3#4}%
        \gdef\cyearmodcHold@{#3#4}}

\def\issn#1{\gdef\theissn{#1}}
\issn{0000-0000}

\def\copyrightline@{\baselineskip2pc
    \rightline{%
        \vbox{\sixrm \textfont2=\sixsy \baselineskip 7pt
            \halign{\hfil##\cr
                \copyright\cyear@\ American Mathematical Society\cr
                 \theissn/\cyearmodc@\ \$1.00 + \$.25 per page\cr }}}}

\def\cyearmodc#1{\gdef\cyearmodc@{\ignorespaces#1\unskip}}

\cvol{000}
\cvolno{0}
\cmonth{}
\cyear{0000}
\cvolyear{0000}
\cpgs{}

%  The copyright block takes the place of the AMS-TeX logo; otherwise the
%  AMSPPT.STY output routine does what we want.

\let\logo@=\copyrightline@

%    Dummy definition of \keyboarder, for now [mjd,1995/04/03]
\def\keyboarder#1{}%

\def\date#1\enddate{\gdef\thedate@{%
\eightpoint Received by the editors \ignorespaces#1\unskip.}}
\def\thedate@{}

\def\dedicatory #1\enddedicatory{\gdef\thededicatory@{{\vskip1pc
  \eightpoint\it \raggedcenter@#1\endgraf}}}
\def\thededicatory@{}

\def\preabstract{}

% presently the "communicated by" line appears only in PROC
% so it is not called in the \endtopmatter definition
% \endtopmatter is redefined in PROC.STY
\def\commby #1\endcommby{\gdef\thecommby@{{\vskip1pc
  \eightpoint \raggedcenter@(Communicated by #1)\endgraf}}}

\let\thecommby@=\empty@

\def\title{\let\savedef@\title
 \def\title##1\endtitle{\let\title\savedef@
   \global\setbox\titlebox@\vtop{\tenpoint\bf
   \raggedcenter@
   \frills@\uppercasetext@{##1}\endgraf}%
 \ifmonograph@ \edef\next{\the\leftheadtoks}%
    \ifx\next\empty@
    \leftheadtext{##1}\fi
 \fi
 \edef\next{\the\rightheadtoks}\ifx\next\empty@ \rightheadtext{##1}\fi
 }%
 \nofrillscheck\title}

\def\author#1\endauthor{\global\setbox\authorbox@
 \vbox{\eightpoint\raggedcenter@
 \uppercasetext@{#1}\endgraf}\relaxnext@ \edef\next{\the\leftheadtoks}%
 \ifx\next\empty@\leftheadtext{#1}\fi}

\def\specialheadfont@{\tenpoint\rm}

\let\specialhead\relax
\outer\def\specialhead#1\endspecialhead{%
  \add@missing\endroster \add@missing\enddefinition
  \add@missing\enddemo \add@missing\endexample
  \add@missing\endproclaim
  \penaltyandskip@{-200}\aboveheadskip
  \begingroup\interlinepenalty\@M\rightskip\z@ plus\hsize
  \let\\\linebreak
  \specialheadfont@\raggedcenter@\noindent
\uppercase{#1}\endgraf\endgroup\nobreak\vskip\belowheadskip}

\let\subsubhead\relax
\outer\def\subsubhead{%
  \add@missing\endroster \add@missing\enddefinition
  \add@missing\enddemo
  \add@missing\endexample \add@missing\endproclaim
  \let\savedef@\subsubhead \let\subsubhead\relax
  \def\subsubhead##1\endsubsubhead{\restoredef@\subsubhead
    \penaltyandskip@{-50}\subsubheadskip
    \varindent@{\def\usualspace{\/{\it\enspace}}%
    \subsubheadfont@##1\unskip\frills@{.\enspace}}}%
  \nofrillscheck\subsubhead}

\long\def\ext #1\endext{\block #1\endblock}

\def\qed{\ifhmode\unskip\nobreak\fi\hfill
  \ifmmode\square\else$\m@th\square$\fi}

\def\xcaheadfont@{\bf}
\outer\def\xca{\let\savedef@\xca \let\xca\relax
  \add@missing\endproclaim \add@missing\endroster
  \add@missing\endxca \envir@stack\endxca
   \def\xca##1{\restoredef@\xca
     \penaltyandskip@{-100}\medskipamount
        \bgroup{\def\usualspace{{\xcaheadfont@\enspace}}%
        \varindent@\xcaheadfont@\ignorespaces##1\unskip
        \frills@{.\xcaheadfont@\enspace}}%
        \ignorespaces}%
  \nofrillscheck\xca}
\def\endxca{\egroup\revert@envir\endxca
  \par\medskip}

\newdimen\rosteritemsep
\rosteritemsep=.5pc

\newdimen\rosteritemitemitemwd
\newdimen\rosteritemitemwd

\newbox\setwdbox
\setbox\setwdbox\hbox{0.}\rosteritemwd=\wd\setwdbox
  \setbox\setwdbox\hbox{0.\hskip.5pc(c)}\rosteritemitemwd=\wd\setwdbox
  \setbox\setwdbox\hbox{0.\hskip.5pc(c)\hskip.5pc(ii)}%
\rosteritemitemitemwd=\wd\setwdbox

\def\roster{%
  \envir@stack\endroster
 \edef\leftskip@{\leftskip\the\leftskip}%
 \relaxnext@
 \rostercount@\z@% Initialize \rostercount@ to 0.
 \def\item{\FN@\rosteritem@}%
 \def\itemitem{\FN@\rosteritemitem@}%
 \def\itemitemitem{\FN@\rosteritemitemitem@}%
 \DN@{\ifx\next\runinitem\let\next@\nextii@\else
  \let\next@\nextiii@\fi\next@}%
 \DNii@\runinitem% If \runinitem occurs, \nextii@ must kill it off.
  {\unskip% This unskips any space before the original \roster.
   \DN@{\ifx\next[\let\next@\nextii@\else
    \ifx\next"\let\next@\nextiii@\else\let\next@\nextiv@\fi\fi\next@}%
   \DNii@[####1]{\rostercount@####1\relax
    \therosteritem{\number\rostercount@}~\ignorespaces}%
   \def\nextiii@"####1"{{\rm####1}~\ignorespaces}%
   \def\nextiv@{\therosteritem1\rostercount@\@ne~}%
   \par@\firstitem@false% Before doing any of this we still change
   \FN@\next@}%      End of definition of \nextii@\runinitem.
 \def\nextiii@{\par\par@% End the present paragraph, change \everypar
  \penalty\@m\vskip-\parskip
  \firstitem@true}%
 \FN@\next@}

\def\rosteritem@{\iffirstitem@\firstitem@false
  \else\par\vskip-\parskip\fi
 \leftskip\rosteritemwd \advance\leftskip\normalparindent
 \advance\leftskip.5pc \noindent
 \DNii@[##1]{\rostercount@##1\relax\itembox@}%
 \def\nextiii@"##1"{\def\therosteritem@{\rm##1}\itembox@}%
 \def\nextiv@{\advance\rostercount@\@ne\itembox@}%
 \def\therosteritem@{\therosteritem{\number\rostercount@}}%
 \ifx\next[\let\next@\nextii@\else\ifx\next"\let\next@\nextiii@\else
  \let\next@\nextiv@\fi\fi\next@}

\def\itembox@{\llap{\hbox to\rosteritemwd{\hss
  \kern\z@ % kern to thwart \unskip in \rom
  \therosteritem@}\hskip.5pc}\ignorespaces}

\def\therosteritem#1{\rom{\ignorespaces#1.\unskip}}

\def\rosteritemitem@{\iffirstitem@\firstitem@false
  \else\par\vskip-\parskip\fi
 \leftskip\rosteritemitemwd \advance\leftskip\normalparindent
 \advance\leftskip.5pc \noindent
 \DNii@[##1]{\rostercount@##1\relax\itemitembox@}%
 \def\nextiii@"##1"{\def\therosteritemitem@{\rm##1}\itemitembox@}%
 \def\nextiv@{\advance\rostercount@\@ne\itemitembox@}%
 \def\therosteritemitem@{\therosteritemitem{\number\rostercount@}}%
 \ifx\next[\let\next@\nextii@\else\ifx\next"\let\next@\nextiii@\else
  \let\next@\nextiv@\fi\fi\next@}

\def\itemitembox@{\llap{\hbox to\rosteritemitemwd{\hss
  \kern\z@ % kern to thwart \unskip in \rom
  \therosteritemitem@}\hskip.5pc}\ignorespaces}

\def\therosteritemitem#1{\rom{(\ignorespaces#1\unskip)}}

\def\rosteritemitemitem@{\iffirstitem@\firstitem@false
  \else\par\vskip-\parskip\fi
 \leftskip\rosteritemitemitemwd \advance\leftskip\normalparindent
 \advance\leftskip.5pc \noindent
 \DNii@[##1]{\rostercount@##1\relax\itemitemitembox@}%
 \def\nextiii@"##1"{\def\therosteritemitemitem@{\rm##1}\itemitemitembox@}%
 \def\nextiv@{\advance\rostercount@\@ne\itemitemitembox@}%
 \def\therosteritemitemitem@{\therosteritemitemitem{\number\rostercount@}}%
 \ifx\next[\let\next@\nextii@\else\ifx\next"\let\next@\nextiii@\else
  \let\next@\nextiv@\fi\fi\next@}

\def\itemitemitembox@{\llap{\hbox to\rosteritemitemitemwd{\hss
  \kern\z@ % kern to thwart \unskip in \rom
  \therosteritemitemitem@}\hskip.5pc}\ignorespaces}

\def\therosteritemitemitem#1{\rom{(\ignorespaces#1\unskip)}}

\def\endroster{\relaxnext@
 \revert@envir\endroster % restore \envir@end
 \par\leftskip@% End the paragraph, and restore the \leftskip.
 \penalty-50% \vskip-\parskip% Add a good break and
 \DN@{\ifx\next\Runinitem\let\next@\relax
  \else\nextRunin@false\let\item\plainitem@% Otherwise, set
   \ifx\next\par% moreover, if \endroster is followed by \par,
    \DN@\par{\everypar\expandafter{\the\everypartoks@}}%
   \else% but if the \endroster isn't followed by a new paragraph,
    \DN@{\noindent\everypar\expandafter{\the\everypartoks@}}%
  \fi\fi\next@}%
 \FN@\next@}

\def\headlinefont@{\sevenpoint}

\def\foliofont@{\sevenrm}

\def\makefootline{\baselineskip18\p@\line{\hss\sevenpoint\folio\hss}}

\def\output@{\shipout\vbox{%
 \iffirstpage@ \global\firstpage@false
\makeheadline
       \pagebody \logo@ \makefootline
\gdef\rightheadline{\hfill \expandafter\iffalse\botmark\fi
  \hfill \llap{\folio}}%
\global\advance\vsize by \topskip % 10pt
     \ifrunheads@
\global\advance\vsize by 18\p@
     \fi
 \else
\ifrunheads@ \makeheadline \pagebody
       \else \pagebody \makefootline \fi
 \fi}%
 \advancepageno \ifnum\outputpenalty>-\@MM\else\dosupereject\fi}

\outer\def\endtopmatter{\add@missing\endabstract
 \edef\next{\the\leftheadtoks}\ifx\next\empty@
  \expandafter\leftheadtext\expandafter{\the\rightheadtoks}\fi
   \ifx\thedate@\empty@\else \makefootnote@{}{\thedate@}\fi
   \ifx\thesubjclass@\empty@\else \makefootnote@{}{\thesubjclass@}\fi
   \ifx\thekeywords@\empty@\else \makefootnote@{}{\thekeywords@}\fi
   \ifx\thethanks@\empty@\else \makefootnote@{}{\thethanks@}\fi
  \jourlogo
  \pretitle
  \begingroup % to localize variant topskip
  \topskip44pt
  \box\titlebox@
  \endgroup
  \preauthor
  \ifvoid\authorbox@\else \vskip2pc\unvbox\authorbox@\fi
  \ifx\thededicatory@\empty@\else
    \thededicatory@\fi
  \preabstract
  \ifvoid\abstractbox@\else
       \vskip1.5pcplus.5pc\unvbox\abstractbox@ \fi
  \ifvoid\tocbox@\else\vskip1.5pcplus.5pc\unvbox\tocbox@\fi
  \prepaper
  \vskip2pcplus1pc\relax
}

%  BOOK REVIEWS

\newcount\affilcount@
\newcount\revrcount@
\newcount\mlcount@
\newcount\revrnum@

{\let\endtopmatter\relax
\let\enddocument=\relax
\gdef\bookrevhook{}% additional macros input by specific style
\gdef\bookrev{%
\gdef\bookrev{\vskip24pt plus 12pt}%

\affilcount@=\z@
\revrcount@=\z@
\mlcount@=\z@

\long\gdef\revtop ##1\endrevtop{\def\revtop@{{\interlinepenalty10000
\tenpoint
\everypar{\hangindent=\parindent}
                \noindent ##1\endgraf}}}%

\gdef\reviewer ##1\endreviewer{\advance\revrcount@\@ne
 \expandafter\gdef\csname revr\number\revrcount@\endcsname{{%
\vbox to 12pt{}%
\leftskip0pt plus1fil
  \rightskip0pt \parfillskip0pt \tenpoint\smc ##1\endgraf}}}%

\gdef\affil ##1\endaffil{%\advance\affilcount@\@ne
 \expandafter\gdef\csname affil\number\revrcount@\endcsname{{%
\leftskip0pt plus1fil
  \rightskip0pt \parfillskip0pt \tenpoint ##1\endgraf}}}%

\gdef\ml ##1\endml{%\advance\mlcount@\@ne
 \expandafter\gdef\csname ml\number\revrcount@\endcsname{{%
\leftskip0pt plus1fil
  \rightskip0pt \parfillskip0pt \tenpoint{\it E-mail
address\,}:\enspace {\tt ##1}\endgraf}}}%

      \gdef\translator@{}%initialize
\gdef\translator ##1\endtranslator{\def\translator@{%
\vskip8pt{\leftskip0pt plus1fil
  \rightskip0pt \parfillskip0pt
  \tenpoint {\rm Translated by\/} ##1\endgraf}}}%

      \gdef\transaffil@{}% initialize
\gdef\transaffil ##1\endtransaffil{\def\transaffil@{%
{\leftskip0pt plus1fil
  \rightskip0pt \parfillskip0pt \tenpoint ##1\endgraf}}}%

     \gdef\transml@{}% initialize
\gdef\transml ##1\endtransml{\def\transml@{{\leftskip0pt plus1fil
  \rightskip0pt \parfillskip0pt \tenpoint{\it E-mail
address\,}:\enspace ##1\endgraf}}}%

  \def\bookrevtext{BOOK REVIEWS}%

      \gdef\bookrevheading{\hrule height\z@
                 \vskip-\topskip
                 \vskip 2.5pc
                 \setbox0=\hbox{\fourteenbf B}
                 \advance\baselineskip by \ht0  %
\global\firstpage@true
\def\rightheadline{}%
\def\\{\break}%
                 \vbox{\raggedcenter@ \fourteenbf \bookrevtext}%
                 \baselineskip45pt\hbox{}\baselineskip0pt}

   \gdef\endtopmatter{}%

  \def\logo@{}%

\font\sixsy=cmsy6

  \def\jourlogo{%
\vbox{%
        \sixrm \textfont2=\sixsy \baselineskip7pt
        \parindent0pt \frenchspacing
AMS JOURNAL STYLE\newline
        Volume {\sixbf\cvol@}, Number \cvolno@, \cmonth@\ \cvolyear@\newline
        \copyright\cyear@\ American Mathematical Society\newline
        \theissn/\cyearmodc@\ \$1.00 + \$.25 per page\par
\vss}}

   \gdef\document{%
        \bookrevheading
        \gdef\bookrevheading{\penalty-100}
\def\\{}%
     \leftheadtext{\bookrevtext}%
\rightheadtext{\bookrevtext}%
\nobreak
        \jourlogo
\nobreak
\vskip8pt  plus5pt minus2pt
\nobreak
\revtop@
\nobreak
\vskip12pt plus5pt minus2pt
\tenpoint }%

   \outer\def\enddocument{\par% \par will do a runaway check for \endref
        \add@missing\endRefs
        \add@missing\endroster \add@missing\endproclaim
        \add@missing\enddefinition
        \add@missing\enddemo \add@missing\endremark \add@missing\endexample
\nobreak
\revrnum@=\z@
\loop\ifnum\revrnum@<\revrcount@\advance\revrnum@\@ne
\csname revr\number\revrnum@\endcsname
   \csname affil\number\revrnum@\endcsname
\csname ml\number\revrnum@\endcsname
%the following is needed to clear the registers
\expandafter\gdef\csname revr\number\revrnum@\endcsname{}%
\expandafter\gdef\csname affil\number\revrnum@\endcsname{}%
\expandafter\gdef\csname ml\number\revrnum@\endcsname{}%
\repeat
\translator@
\transaffil@
\transml@
        }%

\def\no{\gdef\no{\makerefbox\no\keybox@\empty}%
  \gdef\keyhook@{\refstyle A}\no}

\def\rc{\gdef\no{\gdef\no{\makerefbox\no\keybox@\empty}%
  \gdef\keyhook@{\refstyle C}\no}}

\bookrevhook
}% End of \bookrev definition
}% End of group in which \enddocument is non-outer

\catcode`\@=13

%\endinput

%end of ams-j.sty

\catcode`\@=11

\def\jourlogo{%
\def\rightheadline{\vtop to 0pt{%
        \sixrm \baselineskip7pt
        \parindent0pt \frenchspacing
TRANSACTIONS OF THE\newline
AMERICAN MATHEMATICAL SOCIETY\newline
        Volume {\sixbf\cvol@}, Number \cvolno@, \cmonth@\ \cvolyear@\par\vss}}}

\issn{0002-9947}

\let\logo@=\copyrightline@

\catcode`\@=13

%\endinput                                 

%-----------------------------------------------------------------------------
% End of tran.sty
%-----------------------------------------------------------------------------

%the preamble
\NoBlackBoxes

\catcode`\@=11
\def\matrice#1{\left(\null\vcenter {\normalbaselines \m@th
\ialign {\hfil $##$\hfil &&\  \hfil $##$\hfil\crcr
\mathstrut \crcr \noalign {\kern -\baselineskip } #1\crcr
\mathstrut \crcr \noalign {\kern -\baselineskip }}}\right)}

\def\matrize#1{\null\vcenter {\normalbaselines \m@th
\ialign {\hfil $##$\hfil &&\quad  \hfil $##$\hfil\crcr
\mathstrut \crcr \noalign {\kern -\baselineskip } #1\crcr
\mathstrut \crcr \noalign {\kern -\baselineskip }}}}

\def\petitematrice#1{\left(\null\vcenter {\normalbaselines \m@th
\ialign {\hfil $##$\hfil &&\thinspace  \hfil $##$\hfil\crcr
\mathstrut \crcr \noalign {\kern -\baselineskip } #1\crcr
\mathstrut \crcr \noalign {\kern -\baselineskip }}}\right)}

\def\casou#1{\left\{\,\vcenter{\normalbaselines\m@th
      \ialign{$##\hfil$&\quad##\hfil\crcr#1\crcr}}\right.}
\catcode`\@=12

\def\proof{{\sl Proof.\enspace }}

\def\Cont{\mathop{\text {Cont}}}
\def\betacirc{\mathop{\beta_{\text {circ}}}}
\def\betadec{\mathop{\beta_{\text {dec}}}}
\def\betavert{\mathop{\beta_{\text {vert}}}}

\hoffset=1.5cm  \voffset=2cm % to be removed for further printing

\topmatter

\title
A Combinatorial Proof of Bass's Evaluations\\
of the Ihara-Selberg Zeta
Function for Graphs
\endtitle

\author 
Dominique Foata and Doron Zeilberger
\endauthor

\rightheadtext{The Ihara-Selberg Zeta Function for Graphs}

\address
D\'epartement de math\'ematique,
Universit\'e Louis Pasteur, 7, rue Ren\'e-Descartes,
F-67084 Strasbourg, France
\endaddress

\email
foata@math.u-strasbg.fr
\endemail

\address
Department of mathematics, Temple
University, Philadelphia, Pennsylvania 19122
\endaddress

\email
zeilberg@math.temple.edu
\endemail


\abstract
We derive combinatorial proofs of the main two evaluations of the
Ihara-Selberg Zeta function associated with a graph. We give three
proofs of the first evaluation all based on the algebra of Lyndon
words. In the third proof it is shown that the first evaluation is an
immediate consequence of Amitsur's identity on the characteristic
polynomial of a sum of matrices. The second evaluation of the
Ihara-Selberg Zeta function is first derived by means of a
sign-changing involution technique. Our second approach makes use
of a short matrix-algebra argument.
\endabstract

\thanks
The second author was supported in part by N.S.F. and the first
author as a consultant of Zeilberger on his grant.
\endthanks

\date
December 1, 1996
\enddate

\dedicatory
This paper is dedicated to Gian-Carlo Rota, on his millionth$_2$'s
birthday.
\enddedicatory


\endtopmatter

\document
\head 
1. Introduction
\endhead

We are pleased to dedicate the present paper to  our 
{\it bon ma\^\i tre\/} Rota who has been a great promoter of
combinatorial methods. Convinced that combinatorics was hidden in
many branches of mathematics (see, e.g., \cite{14}), he has
successfully  persuaded his followers  to unearth its treasures,
study them for their own sake and propose a fruitful symbiosis with
the  mainstream of mathematics.

Rota's pioneering paper \cite{13} made the {\it M\"obius
function}, and hence its associated {\it zeta function}, a central 
unifying concept in combinatorics and elsewhere. The present paper
is devoted to the calculation of a zeta function, not of a partially
ordered set, as it has been so successfully done in the past by
Rota and his disciples (see, e.g., \cite{17}), but of a tree lattice. 

Digging out those combinatorial treasures is not always an easy
task, since very often a language barrier has to be overcome.  One
such example, that we were fortunate to discover, is Bass's
\cite{3} evaluation of the Ihara-Selberg zeta function for a graph. 
Thanks to his superb and very lucid talk (Temple Mathematics
Colloquium, May 1995) we were introduced to the algebraic set-up
of his derivation and led to the core of his paper. Of great
help also have been his transparencies, copies of which he was kind
enough to send us.

In calculating the zeta function of a tree lattice Bass
\cite{4} was led to determine the following invariant for
a finite connected unoriented graph~$G$. To express his result
he first transformed~$G$ into an {\it oriented} graph
by letting each edge whose ends are vertices~$i$ and~$j$
give rise to two {\it oriented edges} going from~$i$ to~$j$
and from~$j$ to~$i$. Let $c_0$ (resp. $2\,c_1$) be the number
of vertices (resp. of oriented edges). He then introduced the
class~${\Cal R}$ of {\it prime}, {\it reduced} cycles
of~$G$, a class that in general is infinite, and formed the 
product
$$ 
\eta(u)=\prod_{\gamma \in {\Cal R}} (1-u^{|\gamma|}),\leqno(1.1)
$$
where $|\gamma|$ denotes the length of the cycle~$\gamma$. The
product $\eta(u)$ is usually called the {\it Ihara-Selberg function}
associated with the graph~$G$.

His main result was to show that the expansion of $\eta(u)$ as an
infinite series is actually a {\it polynomial} in~$u$ giving two
explicit formulas for it, first as the determinant of a matrix of order
$2\,c_1$ that depends on the {\it successiveness} of the {\it edges}
(a notion that will be defined below), namely 
$$\leqalignno{
\eta(u)&=\det (I-u\,T),&(1.2)\cr
\noalign{\hbox{second, as a product}}
\eta(u)&=(1-u^2) ^{c_1-c_0}\det \Delta(u),&(1.3)\cr}
$$
where $\Delta(u)$ is a matrix of order~$c_0$ that depends on
the {\it connectedness} of the {\it vertices}. The definitions of
$T$ and $\Delta(u)$ will be given in full details later on.

To prove (1.2) and (1.3) Bass makes use of keen algebraic
techniques. In particular, the Jacobi formula
$\det \exp A=\exp {\text {tr}}\,A$ plays a key role in the
derivation of (1.2). As this classical formula has been derived by
combinatorial methods (\cite{6}, \cite{20}), it was
challenging to use those methods to find combinatorial proofs of both
formulas (1.2) and (1.3). This is the purpose of the paper.

With the present combinatorial approach we can show that (1.2) can 
be derived in a more general context. Instead of counting the
cycles by the {\it counter} $u^{|\gamma|}$, we can keep track
of the successiveness property within each cycle~$\gamma$ by
mapping $\gamma$ onto a monomial $\beta(\gamma)$ in the
so-called successiveness variables $b(e,e')$. As we shall see, 
the determinantal expression (1.2) can be derived in three different 
manners, all based on the the algebra of {\it Lyndon words}.

The concept of Lyndon word has been crucial in the
foundations of Free Differential Calculus, initiated by Chen,
Fox and Lyndon \cite{5} and pursued by Sch\"utzenberger \cite{15},
\cite{16} and Viennot \cite{19}. The standard material on the subject
can be found in the book by Lothaire (\cite{8}, chap.~5). Hereafter we
just recall a few basic properties 

Start with a finite nonempty set $X$ supposed to be 
totally ordered and consider the free monoid $X^*$ generated by $X$.
Let $<$ be the lexicographic order on $X^*$ derived from the total
order on $X$. A {\it Lyndon word} is defined to be a nonempty word in
$X^*$ which is prime, i.e., not the power ${l'}^r$  of any other word
$l'$ for any $r\ge 2$, and which is also minimal in the class of its cyclic
rearrangements. Let $L$ denote the set of all Lyndon words. The
following property, due to Lyndon, can be found in \cite{8}, p.~67
(Theorem 5.1.5):

\medskip
\noindent
(1.4) {\sl Each nonempty word $w\in X^*$ can be uniquely written 
as a nonincreasing juxtaposition product of Lyndon words:
$$
w=l_1l_2\ldots l_n,\qquad
l_k\in L,\qquad
l_1\ge l_2\ge \cdots \ge l_n.
$$
}

With each Lyndon word $l$ let us associate a variable denoted by
$[l]$. Assume that all those variables $[l]$ are distinct 
and commute with each other. Furthermore, 
let $\Cal B$ be a square matrix whose entries $b(x,x')$
$(x,x'\in X$) form another set of commuting variables.

If $w=x_1x_2\ldots x_m$ is a nonempty word in $X^*$, define
$$
\betacirc(w):=b(x_1,x_2)b(x_2,x_3)\ldots b(x_{m-1},x_m)b(x_m,x_1)
$$
and $\betacirc(w)=1$ if $w$ is the empty word. Notice that all the 
words in the same cyclic class have the same $\betacirc$-image.
Also  define
$$
\beta([l]):=\betacirc(l)\leqno(1.5)
$$
for each Lyndon word $l$. Now form the $\bold Z$-algebras of formal 
power series in the variables $[l]$ and in the variables 
$b(x,x')$, and by linearity make $\beta$ to be a continuous
homomorphism. It makes sense to consider the product
$$
\Lambda:=\prod_{l\in L} \bigl(1-[l]\bigr)\leqno(1.6)
$$
as well as its inverse $\Lambda^{-1}$. We can also consider the 
images of $\Lambda$ and $\Lambda^{-1}$ under $\beta$. We have
$$\displaylines{
\beta(\Lambda)=\prod_{l\in L} \bigl(1-\beta([l])\bigr);\cr
\noalign{\hbox{and}}
\beta(\Lambda^{-1})=\bigl(\beta(\Lambda)\bigr)^{-1}.\cr
}
$$
We further define two maps $\betadec$ and $\betavert$ 
(``dec" for ``decreasing" and ``vert" for ``vertical") as follows. If
$(l_1,l_2,\ldots, l_n)$ is the nonincreasing factorization of a word $w$
in Lyndon words, as defined in (1.4, let
$$
\betadec(w):=\betacirc(l_1)\betacirc(l_2)\ldots
\betacirc(l_n).
$$
 Now when the $m$ letters of a word $w=x_1x_2\ldots x_m$ are 
rearranged in nondecreasing order, we obtain a word
$\widetilde w=\tilde x_1 \tilde x_2 \ldots \tilde x_m$ called the 
{\it nondecreasing rearrangement} of $w$. Then define
$$
\betavert(w):=b(\tilde x_1,x_1)b(\tilde x_2,x_2)
\ldots b(\tilde x_m,x_m).
$$
Also define $\betadec(w)=\betavert(w):=1$ when $w$ is the empty word.

By convention let $X^*$ denote the {\it sum} of all the words
$w$ ($w\in X^*$) and use the notation
$$
\betadec(X^*):=\sum_{w\in X^*} \betadec(w)
$$
with an analogous notation for $\betavert(X^*)$.

\proclaim{Theorem 1.1} We have the identities
$$\leqalignno{
\beta(\Lambda^{-1})&=\betadec(X^*);&(1.7)\cr
\betadec(X^*)&=\betavert(X^*);&(1.8)\cr
\betavert(X^*)&=\bigl(\det(I-{\Cal B})\bigr)^{-1};&(1.9)\cr
\beta(\Lambda)&=\det(I-{\Cal B}).&(1.10)\cr
}
$$
\endproclaim

Notice that the conjunction of (1.7), (1.8) and (1.9) implies the
identity
$$
\beta(\Lambda^{-1})
=\bigl(\det(I-{\Cal B})\bigr)^{-1}\leqno(1.11)
$$
and therefore (1.10). The proofs of (1.7), (1.8) and (1.9) are given in 
section~2. As we shall see, they are all classical, or preexist in
other contexts. The proof of (1.10) itself is given in section~4. 
Thus we already have two independent proofs of (1.10). 

The direct proof of (1.10) heavily relies on the techniques developed
(or not yet developed) in the algebra of Lyndon words. Section~3  is
then devoted to recalling classical results on  Lyndon words and
proving new~ones. Section~4 contains the construction of an
involution of
$X^*$ that shows that 
$\beta(\Lambda)$ reduces to a {\it finite} sum
$\beta({\bold G})$ that is easily expressible as $\det(I-{\Cal B})$.

Our third proof of (1.10) was suggested to us by Jouanolou \cite{7}
after the first author had discussed the contents of a first version
of the paper during the October 1996 session of the {\it S\'eminaire
Lotharingien}. It is based on a specialization of
Amitsur's identity \cite{2} on the characteristic polynomial of a
finite sum of matrices $A_1+\cdots+A_k$. For each Lyndon word
$l=i_1i_2\ldots i_p$ whose letters belong to the set
$[k]=\{1,2,\ldots,k\}$ let $A_l$ be the matrix product
$A_l:=A_{i_1}A_{i_2}\ldots A_{i_p}$.

Then Amitsur's identity can be stated as
$$
\det(I-(A_1+\cdots +A_k))
=\prod_{l\in L} \det(I-A_l),\leqno(1.12)
$$
where the product is extended over all Lyndon words in the alphabet
$[k]$.

In section 5 we reproduce the (short) proof of Amitsur's identity
(1.12) due to Reutenauer and Sch\"utzenberger \cite {12}. As will be
seen, (1.10) is a mere consequence of (1.12). Thus the
shortest proof of identity (1.10) has
to be borrowed from Classical Matrix Algebra.

In section~6 we show how identity (1.2) fits into
the present context.  It is shown that when $X$ is taken as the
set~$E$ of all oriented edges of the graph~$G$ and each variable
$b(x,x')$ is equal to~0 when the edge~$x'$ is the reverse of~$x$ or
is not the successor of~$x$, and equal to~$u$ otherwise, identity
(1.10) reduces to (1.2).

As named by Bass \cite{3}, the inverse of $\eta(u)$, as given in (1.2), 
is the {\it zeta function} of the underlying tree lattice, so that
$\eta(u)$ itself may be called the {\it M\"obius function} of the tree
lattice. Accordingly, when proving (1.9) (resp. (1.10)) we calculate the
zeta function (resp. the M\"obius function) of the tree lattice.

There is {\it a priori} no extension of (1.3) in which the
information on the edge successiveness can be kept other than a
simple counting of the {\it reduced} prime cycles. We are then left to
prove (1.3) itself, but we present two new proofs, one purely
combinatorial derived in section~7, which is based on the
constructions  of several involutions on words. The second one is
of matrix-algebra nature.

After submitting the present paper for publication in the Fall 1996
our attention was drawn by Ahumada (Mulhouse), who himself published
an early paper on the subject \cite{1}), to the paper by Stark and
Terras that had just appeared~\cite{18}. The latter authors also have
a proof of identity (1.10) when $L$ is restricted to the set of {\it
reduced} prime cycles. Finally, Stanton (Minneapolis) was kind enough
to send us a preprint by Northshield \cite{10} who also has
elementary proofs of both identities (1.2) and (1.3).

\head
2. The zeta function approach
\endhead
When the infinite product $\Lambda$ is developed as an infinite series 
in the variables $[l]$, we get the sum of all the {\it commuting}
monomials $[l_1]\,[l_2]\,\ldots\,[l_n]$, or, equivalently, the sum of the 
nonincreasing words 
$[l_{i_1}]\,[l_{i_2}]\,\ldots\, [l_{i_n}]$
$(l_{i_1}\ge l_{i_2}\ge \cdots\ge l_{i_n})$. Hence, as
$$
\beta(\Lambda^{-1})
=\sum \beta(l_{i_1})\,\beta(l_{i_2})\, \ldots\beta(
l_{i_n})=\sum_{w\in X^*}\betadec(w)=\betadec(X^*),
$$
because of Lyndon's theorem (1.4) and by definition of
$\betadec$, we obtain (1.7).

\medskip
Let $|w|$ denote the length of each word $w\in X^*$.
As both mappings $\betadec$ and $\betavert$ transform a word
{\it of length}~$m$ into a monomial in the variables $b(x,x')$ {\it
of degree}~$m$, identity (1.7) is equivalent to
$$
\sum_{|w|=m}\betadec(w)=\sum_{|w|=m}\betavert(x)
$$
for all $m\ge 0$. Therefore (1.7) is proved if and only if the
following proposition holds.

\medskip
\noindent
(2.1) {\sl There exists a bijection $\Phi$ of $X^*$ onto itself
thaving the following property: if $w=x_1x_2\ldots x_m$ belongs
to~$X^*$, then $\Phi(w)=w'=x'_1x'_2\ldots x'_m$ is a
rearrangement of~$w$ and $\betavert(w')=\betadec(w)$.
}

\medskip
The construction of such a bijection has been given in \cite{4}
(theorem~4.11) and also in \cite{8} p.~198-199. However the
construction must be slightly modified to fit in the present
derivation. We illustrate the construction of the bijection with an
example. Let $X=\{1,2,\ldots,5\,\}$ and
$w=3,4,5,1,2,4,2,1,2,3,1,2,4,2$. The factorization
$(l_1,l_2,\ldots,l_n)$ of~$w$ as a nonincreasing sequence of
Lyndon words (as defined in (1.5)) is
$(3,4,5;\ 1,2,4,2;\ 1,2,3,1,2,4,2)$. For the construction of the
bijection another factorization is used, the {\it decreasing
factorization} $(d_1,d_2,\ldots, d_r)$ of~$w$ simply defined by
cutting~$w$ before every letter~$x$ of~$w$ which is {\it smaller
than or equal to} each letter to its left. With the working example
$(d_1,d_2,\ldots, d_r)=(3,4,5;\ 1,2,4,2;\ 1,2,3;\ 1,2,4,2)$. Notice that
each Lyndon word $l_i$ is the
juxtaposition product of contiguous factors~$d_j$. Moreover
$$\leqalignno{
\betadec(w)&=\betacirc(l_1)\betacirc(l_2)\ldots
\betacirc(l_n)&(2.2)\cr
& =\betacirc(d_1)\betacirc(d_2)\ldots
\betacirc(d_r).\cr}
$$
To obtain $w'$ we form the product of the so-called {\it
dominated circuits} (see \cite{8}, chap.~10)
$$
\Delta(w)
=\left(\matrize{4\ 5\ 3\cr
3\ 4\ 5\cr}\right|
\matrize{2\ 4\ 2\ 1\cr 1\ 2\ 4\ 2\cr}
\left|\matrize{2\ 3\ 1\cr 1\ 2\ 3\cr}\right|
\left.\matrize{2\ 4\ 2\ 1\cr
1\ 2\ 4\ 2\cr}\right).
$$
In $\Delta(w)$ the top word in each factor is obtained from the
bottom factor~$d_j$ by making a right to left cyclic shift
of~$d_j$.

Next we reshuffle the columns of $\Delta(w)$ in such a way that
the mutual order of two columns with the same top entry is not
modified but the top row becomes nonincreasing:
$$
\left(\matrize{5\ 4\ 4\ 4\ 3\ 3\ 2\ 2\ 2\ 2\ 2\ 1\ 1\ 1\cr
4\ 3\ 2\ 2\ 5\ 2\ 1\ 4\ 1\ 1\ 4\ 2\ 3\ 2\cr}\right).
$$
The resulting bottom word is the word $\Gamma^{-1}(\Delta(w))$
as described in \cite{8}, p. 199, except the construction has been
given with the {\it reverse} order of~$X$.

Now exchange top and bottom words and rewrite the resulting biword
from right to left:
$$
\left(\matrize{2\ 3\ 2\ 4\ 1\ 1\ 4\ 1\ 2\ 5\ 2\ 2\ 3\ 4\cr
1\ 1\ 1\ 2\ 2\ 2\ 2\ 2\ 3\ 3\ 4\ 4\ 4\ 5\cr}\right).
$$
Finally, reshuffle the columns of the last biword so that the top
word becomes {\it nondecreasing}, still keeping the mutual order
of any two columns having the same top entry invariant:
$$\displaylines{\matrize{\cr w'={}\cr}
\left(\matrize{1\ 1\ 1\ 2\ 2\ 2\ 2\ 2\ 3\ 3\ 4\ 4\ 4\ 5\cr
2\ 2\ 2\ 1\ 1\ 3\ 4\ 4\ 1\ 4\ 2\ 2\ 5\ 3\cr}\right).\cr
\noalign{\hbox{Then $w'$ is defined to be the bottom word of the
above biword. Moreover}}
\betadec(w)=\betavert(w').\cr
\noalign{\hbox{With the working example the latter monomial is
equal to}}
b(1,2)^3\,b(2,1)^2\,b(2,3)\,b(2,4)^2\,b(3,1)\,
b(3,4)\,b(4,2)^2\,b(4,5)\,b(5,3).\cr}
$$

Identity (1.9) is essentially the MacMahon Master Theorem identity
(see \cite{9}, p. 93-96, or \cite{4}, chap.~5). This achieves the
proof of~(1.11). Notice that the combination of (1.7), (1.8) and
(1.10) provides a new proof of the Master Theorem identity.

\head
3. Lyndon and Donlyn words
\endhead
As already defined in the introduction a {\it Lyndon word} is a
nonempty word in $X^*$ which is prime and also minimal in its
class of cyclic rearrangements. Let $L$ denote the set of all
Lyndon words. The following properties (3.1)--(3.3) can be found in
\cite{8}, pp. 65 and 66 (Propositions 5.1.2 and 5.1.3):

\medskip
\noindent
(3.1) {\sl A nonempty word in $X^*$ is a Lyndon word if and only
if it is strictly smaller that any of its proper right factors.}

\smallskip
\noindent
(3.2) {\sl A nonempty word in $X^*$ is a Lyndon word if and only
if it is of length one or the juxtaposition product $lm$ of two
Lyndon words $l$, $m$ such that $l<m$.}

\smallskip
Let $l$ be a Lyndon word; if $|l|\ge 2$ let $m_0$ be the proper
right factor of maximal length such that $m_0\in L$. Write
$l=l_0m_0$. The factorization $(l_0,m_0)$ of~$l$ is called the
{\it standard factorization} of~$l$.

\medskip
\noindent
(3.3) {\sl If $(l_0,m_0)$ is the standard factorization of a
Lyndon word~$l$ of length $|l|\ge 2$, then $l_0$ is also a
Lyndon word and $l_0<l_0m_0<m_0$.}

\medskip
\medskip
We will also need the following two properties, apparently not
stated in the standard texts, but essential in our derivation.

\medskip
\noindent
(3.4) {\sl A factorization $(l_0,m_0)$ of a
Lyndon word~$l$ into two nonempty factors is the standard
factorization of~$l$ if and only if
$m_0l_0$ is the second smallest cyclic rearrangement of $l$ $($the
smallest one being
$l$ itself\/$)$.}

\medskip
For obvious reasons we shall call the word $m_0l_0$ a {\it
Donlyn} word. We reproduce the short proof kindly provided
by Perrin~\cite{11}.

Notice that if $(l_0,m_0)$ is the standard factorization of~$l$,
then $m_0$ is necessarily the smallest proper right factor of~$l$. Let
$(l_1,m_1)$ be another factorization of~$l$. Either $m_1$ does not
start with~$m_0$ and then
$m_0l_0<m_1l_1$, or $m_1=m_0m_2$ for some word~$m_2$. In
the latter case, as $m_2$ is a proper right factor of~$l$, we
have $l<m_2$ and then $m_0l_0<m_0l<m_0m_2l_1=m_1l_1$.
The converse is immediate.\qed

\medskip
\noindent
(3.5) {\sl Let $l$, $m$ be two Lyndon words such that $l<m$.
Then $(l,m)$ is the standard factorization of $lm$ if and only if
$m$ is less than each of the cyclic rearrangements of~$l$ other
than~$l$.}

\medskip
\proof
Assume that $(l,m)$ is the standard factorization of~$lm$ and
let $l=l'l''$ with both $l'$ and $l''$ nonempty. If $l''=mm''$,
then $m<l''l'$. If $l''$ does not start with $m$, then
$m<l''$; otherwise, we would have $l''<m$ and then
$l''m<m$ which contradicts the fact that $m$ is the smallest
proper right factor of~$lm$. Now if $m=l''m'$, then 
$l''m'=m<l''m=l''l''m'$ implies $m'<l''m'=m$ and this
contradicts the fact that $m$ is a Lyndon word. Accordingly,
$m$ cannot start with~$l''$ and the inequality $m<l''$
implies $m<l''l'$.

Conversely, suppose that $m$ is less than each of the 
cyclic rearrangements of~$l$ other than~$l$. If $(l,m)$ is not
the standard factorization of~$lm$, then $l=l'l''$, with
$l'$, $l''$ nonempty, $l''m\in L$ and $l''m<m$. By assumption, we
also have $m<l''l'$. Therefore, $l''m<m<l''l'$. This implies that
$m=l''m'$ with $m<m'<l'$, so that $m<m'<l'<l$. But the
inequality $m<l$ cannot hold as $lm$ is a Lyndon word.\qed

\head
4. The M\"obius function approach
\endhead
The {\it content} of a word~$w$ is defined as the set 
$\Cont (w)$ of all distinct letters occurring in~$w$.
A nonempty word of $X^*$ is said to be {\it multilinear}, if all
its letters are distinct. Two words $w$ and $w'$ are said to be
{\it disjoint}, if they have no letter in common.

Denote by $[L]$ the set of all commuting variables $[l]$
associated with each Lyndon word~$l$. If $w$ is a prime word,
it is the cyclic rearrangement of a unique Lyndon
word~$l$. We will also write $[w]=[l]$, regarding each
variable~$[l]$ as being associated with the {\it class of cyclic
rearrangements} of the word~$l$. We next form the free
Abelian monoid ${\text {Ab}}[L]$ generated by
$[L]$ and consider the following sequence
${\text {Ab}}[L]\supset{\Cal D}\supset {\Cal G}$ defined as follows.
Each monomial $\pi=[l_1]\,[l_2]\,\ldots\,[l_r]$ belongs to
$\Cal D$, if and only if the Lyndon words $l_1$, $l_2$,
\dots~, $l_r$ are all {\it distinct}. It belongs to~$\Cal G$ if
furthermore every element $x\in X$ occurs {\it at most once} in 
the set $\Cont(\pi)=\Cont(l_1l_2\ldots l_r)$. In such a case all
the Lyndon words $l_k$ are necessarily multilinear. As $X$ is
finite, the set~$\Cal G$ is necessarily {\it finite}. Moreover each
element
$\pi\in {\Cal G}$ may be regarded as a {\it permutation} of the set
$\Cont(\pi)\subset X$. The number~$r$ of factors in $\pi$ is called the
{\it degree} of~$\pi$ and denoted by $\deg\pi$.

The expansion of~$\Lambda$ (defined in (1.6)) is the infinite series
$$
\leqalignno{
\Lambda&=\sum_{\pi\in {\Cal D}}(-1)^{\deg\pi}\pi.&(4.1)\cr
\noalign{\hbox{We can also form the {\it polynomial}}}
{\bold G}&:=\sum_{\pi\in {\Cal G}}(-1)^{\deg\pi}\pi.&(4.2)\cr
}
$$
The definition of the homomorphism $\beta$ was given in (1.5). 

\proclaim{Theorem 4.1} We have the identity
$$
\beta(\Lambda)=\beta({\bold G}),\leqno(4.3)
$$
so that $\beta(\Lambda)$ is a polynomial.
\endproclaim

The proof of Theorem 4.1 is based on an {\it involution}
$\pi\mapsto \pi'$ of ${\Cal D}\setminus {\Cal G}$ 
such that $\deg \pi+\deg \pi'=0$ mod~2 that is defined as
follows. 

\medskip
{\it Construction of the involution}. Say that
$\pi=[l_1]\,[l_2]\,\ldots\,[l_r]$ is a {\it good companion} if it
belongs to $\Cal G$. If 
$\pi=[l_1]\,[l_2]\,\ldots\,[l_n]$ is a bad companion (an
element of ${\Cal D}\setminus {\Cal G}$), let $x$ be the
{\it smallest} letter that occurs more than once
in $l_1l_2\ldots l_r$. If 
$l_i$ contains~$x$, let $xu_1$, $xu_2$, \dots~, $xu_s$ be
the list of all cyclic rearrangements of~$l_i$ that
start with~$x$. Write such a list for each of the words 
$l_1$, $l_2$, \dots~, $l_r$ and combine all
those lists. It is essential to notice that all the elements
in the list are {\it distinct}, because it is so for all the cyclic
rearrangements of a Lyndon word and by assumption all the
Lyndon words $l_1$, $l_2$, \dots~, $l_r$ are themselves distinct.

Now choose a total order on~$X$ such that $x=\min X$ and
consider the lexicographic order on $X^*$ with respect to that
total order. Furthermore, write the previous list in {\it
increasing order}
$$
{\text {List}}(\pi)=(xu_1, xu_2, xu_3, \ldots\,)\leqno(4.4)
$$
and consider the {\it smallest two elements} $xu_1$, $xu_2$.
Either they come from the same factor $l_i$ (case (i)), or
from two different factors (case (ii)).

In case (i) write $xu_1=xvxw$, $xu_2=xwxv$ ($v,w\in X^*$) so
that 
$[l_i]=[xvxw]=[xwxv]$. Then define
$$
\pi=[l_1]\,[l_2]\,\ldots [l_r]\,
\mapsto\,
\pi'=[l_1]\,\ldots\,[l_{i-1}]\,[xv]\,[xw]\,
[l_{i+1}]\,\ldots [l_r].
$$

In case (ii) suppose $[xu_1]=[l_1]$, $[xu_2]=[l_2]$.
Then define
$$
\pi=[l_1]\,[l_2]\,\ldots \,[l_r]\,
\mapsto\,
\pi'=[xu_1xu_2]\,
[l_3]\,\ldots\,[l_r].
$$

In case (i) the word $xu_1=xvxw$ that is first
in ${\text {List}}(\pi)$ is necessarily a Lyndon word (with respect to
the latter total order on~$X$). Furthermore, the pair $(xv,xw)$ is
the standard factorization of
$xvxw$ by Property (3.4). Therefore, both $xv$, $xw$ are Lyndon
words and accordingly prime by Property~(3.3).

On the other hand, as $xvxw$ and $xwxv$ are the smallest two
elements in ${\text {List}}(\pi)$ and since
$xv<xvxw<xw<xwxv<xu_3<\cdots$, both $xv$ and $xw$ are
smaller than all the other words $xu_k$ for each $k\ge 3$.
It also follows from Property~(3.4) that $xw$ is less than all
the cyclic rearrangements of~$xv$ other than~$xv$.
Accordingly, the smallest two elements in ${\text {List}}(\pi')$
are $xv$ and $xw$. Consequently, $(\pi')'=\pi$.

\smallskip
In case (ii) the two words $xu_1$ and $xu_2$ coming
from two different factors are necessarily Lyndon words. As
$xu_1<xu_2$, Property (3.2) implies that $xu_1xu_2$ is also a
Lyndon word and therefore is prime. On the other hand, as
$xu_1<xu_1xu_2<xu_2$, the word $xu_1xu_2$ is 
less than all cyclic rearrangements of $xu_1$ 
that may occur in ${\text {List}}(\pi)$ other than all the words
$xu_k$ for each
$k\ge 3$, in particular it is less than all cyclic
rearrangements of $xu_1$ that may occur in ${\text {List}}(\pi)$ other
than $xu_1$. It follows from Property~(3.5) that
$(xu_1,xu_2)$ is the standard factorization of $xu_1xu_2$.
Accordingly, $xu_1xu_2$ and $xu_2xu_1$ are the smallest two
elements in ${\text {List}}(\pi')$. Hence $(\pi')'=\pi$.

This shows that $\pi\mapsto \pi'$ is a well defined involution
of ${\Cal D}\setminus {\Cal G}$. Moreover it satisfies
$$
\beta(\pi)=\beta(\pi').
$$
Therefore (4.3) holds.\qed

\medskip
Let $\pi=[l_1]\,[l_2]\,\ldots\, [l_r]$ be a monomial
belonging to ${\Cal G}$. As noted before, $\pi$ may be
regarded as a {\it permutation} of~$\Cont(\pi)$. The set~$\Cal
G$ is then the set of all permutations of {\it subsets} of~$X$. The
summand $(-1)^{\deg \pi}\beta(\pi)$ in $\beta({\bold G})$ is then the
term in the expansion of $\det(I-{\Cal B})$ associated with
the permutation~$\pi$ (see, e.g., \cite{20}, \S\kern2pt 1). Thus
$$
\beta({\bold G})=\sum_{\pi\in {\Cal G}}
(-1)^{\deg \pi}\,\beta(\pi)=\det (I-{\Cal B}).\leqno(4.5)
$$
This yields identity (1.10) in view of Theorem 4.1.

\head
5. Amitsur's identity
\endhead
Reutenauer and Sch\"utzenberger \cite{12} gave the following short
proof of the Amitsur identity~(1.12): Lyndon's factorization
theorem~(1.4) may be expressed as
$\prod_l (1-l)^{-1}=X^*=(\sum_w w)$ $(w\in X^*)$ where the
product is taken over all Lyndon words in nonincreasing order. As
$X^*=(1-X)^{-1}=(\sum_w w)$ $(w\in X^*)$, we can deduce
$(1-X)^{-1}=\prod_l (1-l)^{-1}$. Now form the inverse of the
latter identity and replace
$X$ by a set of matrices $\{A_1,\ldots, A_k\}$. Taking the
determinant of both sides yields identity~(1.12).

Next Amitsur's identity (1.12) specializes into (1.10) in the following
manner. Let
$N=2c_1$,
$k=N\times N$ and consider the lexicographic order on the pairs $(i,j)$
$(1\le i,j\le N)$. If $(i,j)$ is the $m$-th pair,
let $A_m$ be the matrix whose entries are all null except the
$(i,j)$-entry which is equal to $b(i,j)$. Then $A_1+\cdots
+A_k={\Cal B}$.

Consider a word $l=(i_1,j_1)(i_2,j_2)\ldots (i_p,j_p)$ in the alphabet
$\{(1,1),\ldots, (N,N)\}$. If
$j_1=i_2$, $j_2=i_3$, \dots~, $j_{p-1}=i_p$, then $A_l$ is the
matrix whose all entries are null except the $(i_1,j_p)$-entry which
is equal to $b(i_1,i_2)b(i_2,i_3)\cdots b(i_{p-1},i_p)b(i_p,j_p)$.
If  the above contiguity relations for the entries~$b(i,j)$ do not hold,
$A_l$ is zero. 

Now remember that $\det(I-A_l)$ is the alternating sum of the {\it
diagonal} minors of the matrix~$A_l$. Accordingly, when the word~$l$
satisfies the above contiguity relations {\it and}
$j_p=i_1$, we have
$$
\det(I-A_l)=
1-b(i_1,i_2)b(i_2,i_3)\cdots b(i_{p-1},i_p)b(i_p,i_1).
$$
In the other cases, $\det(I-A_l)=1$.

The infinite product in (1.12) can then be restricted to the Lyndon
words $l=(i_1,j_1)(i_2,j_2)\ldots (i_p,j_p)$ in the alphabet
$\{(1,1),\ldots, (N,N)\}$ satisfying
$j_1=i_2$, $j_2=i_3$, \dots~, $j_{p-1}=i_p$ et $j_p=i_1$.
But those words are in bijection with the Lyndon words
$i_1i_2\ldots i_p$ in the alphabet~$[N]$. This proves identity~(1.10).

\head
6. Bass's results
\endhead
As said in the introduction Bass's calculations deal with an
oriented graph having $c_0$ vertices labelled
$1,2,\ldots,c_0$ and
$2\,c_1$ oriented edges. Notice that each loop around
vertex~$i$ in the original unoriented graph gives
rise to {\it two oriented} loops around~$i$ in the oriented
graph. Each oriented edge~$e$ going from vertex~$i$, called the {\it
origin\/} of~$e$, to vertex~$j$, called the {\it end} of~$e$, has a
unique {\it reverse} edge going from $j$ to $i$ that will be denoted
by~$J(e)$ or by~$\overline e$. Let $V$ be the set of vertices and $E$
be the set of oriented edges so that $\# V=c_0$ and  $\#E=2c_1$.

Say that an oriented edge $e'$ is a {\it successor} of an
oriented edge~$e$, if the end of~$e$ and the origin of~$e'$
coincide. An {\it oriented path} from vertex~$i$ to vertex~$j$
is a linear sequence of oriented edges $e_1e_2\ldots e_m$
$(m\ge 1)$ such that for every $k=1,2,\ldots,m-1$ the edge
$e_{k+1}$ is a sucessor of~$e_k$; moreover the origin of $e_1$
is~$i$ while the end of~$e_m$ is~$j$. The integer~$m$ is the
{\it length} of the oriented path. It will be convenient to consider
the free monoid~$E^*$ generated by the edge set~$E$ and see the
oriented paths as particular elements of~$E^*$. 

When $j=i$ the oriented path is called a {\it pointed cycle}. The
oriented path $e_1e_2\ldots e_m$ is said to be {\it reduced}, if
$J(e_1)\not=e_2$, $J(e_2)\not=e_3$, \dots~,
$J(e_{m-1})\not=e_m$, $J(e_m)\not=e_1$.
A pointed cycle is said to be {\it prime}, if it cannot
be expressed as the product $\delta^r$ of a given
pointed cycle~$\delta$ for any $r\ge 2$. 

Two pointed cycles~$\delta$ and ~$\delta'$ are said to be
(cyclically) {\it equivalent}, if they are cyclic rearrangements
of each other, i.e., if they can be expressed as words
$\delta=e_1e_2\ldots e_m$ and $\delta'=e_ke_{k+1}\ldots
e_me_1\ldots e_{k-1}$ in $E^*$ for some $k$ $(1\le k\le m)$.
Each equivalence class is called a {\it cycle}. The cycle
containing the pointed cycle~$\delta$ will be denoted by
$[\delta]$. This notation will not conflict with our previous
notation for the variables~$[l]$ as we shall see.

If a pointed cycle is prime (resp. reduced, resp.
of length $m$), all the elements in its equivalence class are prime
(resp. reduced, resp. of length~$m$). We can
then speak of {\it prime, reduced cycles}. The {\it length} of a
cycle~$\gamma$ will be denoted by $|\gamma|$. Let ${\Cal P}$
(resp. ${\Cal R}$) denote the set of all {\it
prime} (resp. {\it prime} and {\it reduced\/}) cycles. The
ingredients of (1.1) are then fully defined.

The further notions introduced by Bass are the following.

(i) For each $i=1,2,\ldots,c_0$ let $E_i$ (resp. ${\Cal L}(E_i)$) be
the set of the oriented edges going out of vertex~$i$ (resp. the
vector space spanned by the basis~$E_i$). The {\it outer degree}
of vertex~$i$ is the number of oriented edges going
out of~$i$. Let $Q(i)$ be equal to that outer degree minus one,
so that, as the graph is assumed to be connected, $Q(i)\ge 0$.
Let ${\Cal Q}$ be the diagonal matrix
${\text {diag}}(Q(1),\ldots,Q(c_0))$.
With those notations we have: $\dim {\Cal L}(E_i)=Q(i)+1$. The direct
sum of all the ${\Cal L}(E_i)$'s will be denoted by~${\Cal L}(E)$, so
that
${\text {dim}}\,{\Cal L}(E)=\sum\limits_{i=1}^n (Q(i)+1)=2\,c_1$.

\smallskip
(ii) The {\it successiveness map} ``Succ" is defined as follows:
let $e$ be an oriented edge going from vertex~$i$ to
vertex~$j$. Then
$$ 
{\text {Succ}}(e):=\sum_{e'\in E_j} e'.\leqno(6.1) 
$$
In other words, ${\text {Succ}}(e)$ is the sum of all the
successors of~$e$. The mappings Succ and the reverse map~$J$
may be regarded as endomorphisms of~${\Cal L}(E)$. Then 
$T={\text {Succ}}-J$ is the endomorphism occurring in formula (1.2).

\smallskip
(iii) The {\it connectedness matrix} 
${\Cal K} =(K(i,j))$ $(1\le i,j \le c_0)$. Let $E_{i,j}$ be the set of
all oriented edges going from vertex~$i$ to vertex~$j$. Then
$K(i,j):=\left|E_{i,j}\right|$. Notice that $K(j,i)=K(i,j)$ and
$K(i,i)\ge 2$ if there is a loop around~$i$. 
The matrix $\Delta(u)$
occurring in (1.3) is the matrix
$$
\Delta(u)=I-u\,{\Cal K}+u^2\,{\Cal Q}.\leqno(6.2)
$$

To recover the first evaluation (1.2) we have to take the
following ingredients: 

(i) $X=E$, the set of oriented edges;

(ii) ignore each variable $b(e,e')$ when
$e'$ is {\it not} a successor of~$e$ or when $e'=J(e)$ (mapping it
to~0) and make all the other variables equal to~$u$. Call
$\beta_u$ the corresponding homomorphism~$\beta$.

If a cycle is prime, it contains a unique pointed cycle which is also
a Lyndon word~$l$. We then denote the cycle by $[l]$. We have
$$\displaylines{
\beta_u([l])=\casou{u^{|l|},&if
$[l]$ is reduced;\cr
0,&otherwise;\cr}\cr
\noalign{\hbox{and then}}
\beta_u(\Lambda)=\prod_{\gamma\in{\Cal R}}(1-u^{|\gamma|}).\cr}
$$
Also if $\pi=[l_1]\,[l_2]\,\ldots\,[l_r]$ is a monomial whose
components are prime reduced cycles, we have
$\beta_u(\pi)=u^{|\Cont\pi|}$.
Let $\Cal H$ be
the set of the monomials 
$\pi=[l_1]\,[l_2]\,\ldots\,[l_r]$ such that each $[l_k]$ is a prime
reduced cycle and every edge occurs at most once in
$l_1l_2\ldots l_r$ and let
$$\displaylines{
{\bold H}:=\sum_{\pi\in {\Cal H}}(-1)^{\deg \pi} \pi.\cr
\noalign{\hbox{Then}}
\beta_u({\bold G})=\beta_u({\bold H}).\cr
\noalign{\hbox{so that (4.3) becomes}}
\rlap{(6.3)}\hfill
\beta_u(\Lambda)=\beta_u({\bold H}).\hfill\cr}
$$
As $\det (I-{\Cal B})$ reduces to $\det(I-u\,T)$, formula (4.5)
becomes 
$$
\beta_u({\bold H})=\sum_{\pi\in{\Cal H}}(-1)^{\deg\pi}
u^{|\Cont\pi|}=\det(I-u\,T).\leqno(6.4)
$$
\head
7. A Purely Combinatorial Proof of Formula (1.3)
\endhead
Our purpose is to give a combinatorial proof of the identity
$$
(1-u^2)^{{1 \over 2}
\vert E \vert + \vert V \vert} \beta_u({\bold H})=
(1-u^2)^{\vert E \vert} \det \Delta(u)  ,\leqno(7.1)
$$
which is obviously equivalent to (1.3) because of (6.4). 
The determinant $\Delta(u)$ was defined in (6.2).

Our strategy will be to introduce a class of permutation graphs
with colored edges, called {\it chaps} and consider the sum of the  weights
of all those chaps. That sum will be computed in two
different ways. We will soon define {\it polite chaps} and later
{\it good chaps}. It will be shown that the weighted sum of the
impolite chaps is zero, as well as the weighted sum of the bad
chaps. This is achieved by defining appropriate involutions that
will partition all the impolite chaps into pairs each of whose
members' weight is the negative of the other, and similarly for
the bad chaps. It will then follow that the sum of the weights
of the polite chaps equals the sum of the weights of the
good chaps. The former will turn out to be the right side of (7.1)
while the latter will turn out to be the left side of (7.1).

\subhead{\rm 7.1.}\enspace
Introducing Chaps\endsubhead
Suppose that the set~$E$ of all oriented edges of~$G$ is totally
ordered. A {\it chap} may be seen as a {\it permutation graph}
$\text {Ch}$ (i.e., a collection of disjoint cycles) whose vertices
--- call them {\it supervertices} --- are the vertices and the
edges of the original graph, i.e., the elements of $V\cup E$, and
whose edges --- call them {\it superedges} --- are {\it
colored} in the following sense. Let $e$, $e'$ be two oriented
edges (not necessarily distinct) going out of the same
vertex~$i$, let $j$ be the end of~$e$ and let $e''$ be a
successor of~$e$ (its origin is then vertex~$j$). By definition
the only possible {\it colored superedges} of a chap are the
following 
$$\displaylines{i \buildrel 1\over
\longrightarrow i ;\quad i \buildrel 2\over \longrightarrow e ;\quad
e\buildrel 3a\over \longrightarrow i ;\quad
e \buildrel 3b\over \longrightarrow i ;\cr
e \buildrel 4a\over \longrightarrow j ;\quad
e \buildrel 4b\over \longrightarrow j ;\quad
e \buildrel 5\over \longrightarrow e ;\quad
e \buildrel 6\over \longrightarrow e' ;\quad
e\buildrel 7\over \longrightarrow e''.\cr}
$$
With each of the nine colors is associated a weight as shown in the
next table:
$$\vbox{\halign{\vrule\vrule depth 5pt height 12pt width 0pt\
\hfil$#$\hfil\  &\vrule\ \hfil$#$\hfil\ 
&\vrule\ \hfil$#$\hfil\ &\vrule\ \hfil$#$\hfil\
&\vrule\ \hfil$#$\hfil\  &\vrule\ \hfil$#$\hfil\ &\vrule\ 
\hfil$#$\hfil\ &\vrule\ \hfil$#$\hfil\  &\vrule\ \hfil$#$\hfil\
&\vrule\ \hfil$#$\hfil\ \vrule\cr
\noalign{\hrule}
{\text {Color}}&1&2&3a&3b&4a&4b&5&6&7\cr 
\noalign{\hrule}
{\text {Weight}}&1-u^2&u(1-u^2)&u&-u&1&-1&1-u^2&u^2&-u\cr
\noalign{\hrule}
}}
$$
The {\it weight} of a chap is defined to be the product
of the weights of the superedges times the {\it sign} of the graph
permutation.

A chap is {\it polite} if its superedges are of the form:
$$i \buildrel 1\over \longrightarrow i ;\quad
i \buildrel 2\over \longrightarrow e ;\quad
e \buildrel 3b\over \longrightarrow i ;\quad
e \buildrel 4a\over \longrightarrow j ;\quad
e \buildrel 5\over \longrightarrow e;
$$
where $e$ is an edge of origin~$i$ and end~$j$. 
A chap that is not polite will be called {\it impolite}.
If a chap is impolite, there exists a vertex~$i$ such that one of the
following conditions holds from some edges~$e$, $e'$, $e''$:
(A) $e \buildrel 3a\over \longrightarrow i$, $e\in E_i$;
(B) $e \buildrel 4b\over \longrightarrow j$, $e\in E_i$;
(C) $e \buildrel 6\over \longrightarrow e'$, $e\in E_i$;
(D) $e \buildrel 7\over \longrightarrow e''$, $e\in E_i$.
Denote by $i$ the smallest such a vertex and 
let $e$ be the {\it smallest} oriented edge
in~$E_i$ which is the origin of a superedge colored $3a$, $4b$, 6
or~7. Accordingly, one the following six conditions holds:

\smallskip
(1) $e \buildrel 3a\over \longrightarrow i \buildrel 2\over
\longrightarrow e'$; 

(2) $e \buildrel 4b\over \longrightarrow 
j \buildrel 2\over \longrightarrow e'$;

(3) $e \buildrel 6\over \longrightarrow e'$, 
$i \buildrel 1\over \longrightarrow i$; 

$(3')$ $e \buildrel 6\over \longrightarrow 
e'$, $e''\buildrel x\over \longrightarrow i
\buildrel 2\over \longrightarrow e'''$ 
with $x=3a,\,3b,\,4a$ or~$4b$.

(4) $e \buildrel 7\over \longrightarrow e'$, 
$j \buildrel 1\over \longrightarrow j$, $e'\in E_j\,$; 

$(4')$ $e \buildrel 7\over \longrightarrow e'$, 
$e''\buildrel x\over \longrightarrow j
\buildrel 2\over \longrightarrow e'''$, $e'\in E_j$
with $x=3a,\,3b,\,4a$ or~$4b$.

\smallskip
\noindent
If (1) (resp. (2)) occurs within an impolite chap~Ch,
transform~Ch into another (impolite) chap~${\text {Ch}}'$ by
replacing occurrence~(1) (resp. (2)) by occurrence~(3) (resp. (4))
and conversely. Finally, if $(3')$ occurs, perform the change:
$e \buildrel 6\over \longrightarrow e'''$, 
$e''\buildrel x\over \longrightarrow i
\buildrel 2\over \longrightarrow e'$ and if $(4')$ occurs,
perform the change
$e \buildrel 7\over \longrightarrow e'''$ and
$e''\buildrel x\over \longrightarrow j
\buildrel 2\over \longrightarrow e'$. Those
changes preserve the absolute value of the weight and reverse
its sign.

It follows that the sum of the weights of all impolite chaps
is zero. Hence the sum of the weights of all chaps equals the
sum of the weights of the polite chaps. We will now proceed
to compute it. 

\subhead{\rm 7.2.}\enspace
The sum of the weights of the polite chaps\endsubhead
Each polite chap consists of cycles 
where superedges 2 and 4a intertwine
$$
i_1\buildrel 2\over \longrightarrow
e_1\buildrel 4a\over \longrightarrow
i_2\buildrel 2\over \longrightarrow
e_2\buildrel 4a\over \longrightarrow i_3\ 
\cdots\ 
i_k\buildrel 2\over \longrightarrow
e_k\buildrel 4a\over \longrightarrow
i_1
$$
as well as $2$-cycles of the form $i\buildrel 2\over
\longrightarrow e \buildrel 3b\over \longrightarrow i$, the
other vertices and edges being fixed points:
$i\buildrel 1\over \longrightarrow i$,
$e\buildrel 5\over \longrightarrow e$.

A cycle of the first kind has weight $u^k(1-u^2)^k$, while a
cycle of the second kind has weight $-u^2(1-u^2)$. To the
product of all these cycles we must multiply by $(1-u^2)$ raised
to the power of the number of remaining edges and vertices. 
Let $V_1$ (resp. $V_2$, resp. $V_3$) be the set of vertices
belonging to the cycles of the first kind (resp. of the second
kind, resp. of the form $e\buildrel 5\over \longrightarrow e$).
Since each cycle of the first kind has the same number of
vertices and edges, and each cycle of the second kind has one
vertex and one edge, the total weight is
$$
(1-u^2)^{\vert E \vert}\times u^{\left|V_1\right|}
\times (-u^2)^{\left|V_2\right|}\times
(1-u^2)^{\left|V_3\right|}.
$$
This is the same
as $(1-u^2)^{\vert E \vert+ \vert V \vert }\times 
(u/(1-u^2))^{\left|V_1\right|}\times
(-u^2/(1-u^2))^{\left|V_2\right|}$.

\goodbreak
Remember that $E_{i,j}$ denote the set of all oriented edges in
the graph~$G$ going from vertex~$i$ to vertex~$j$ and
$|E_{i,j}|=K(i,j)$. A polite chap~Ch is then characterized by a
sequence
$(V_1,V_2,V_3,\sigma,f,g)$, where

(i) $(V_1,V_2,V_3)$ is a partition of the vertex set $V$ in
disjoint subsets;

(ii) $\sigma$ is a permutation of $V_1$;

(iii) $f:V_1\rightarrow E$ is a mapping such that
$[\sigma(i)=j]\Rightarrow [f(i)\in E_{i,j}]$;

(iv) $g:V_2\rightarrow E$ is a mapping such that $g(i)\in E_i$.

\noindent
Write $\alpha=u/(1-u^2)$ and $\beta=-u^2/(1-u^2)$. As the
sign of~$\pi$ is given by
$\varepsilon(\sigma)\,(-1)^{|V_1|+|V_2|}$, the sum of the
weights of the polite chaps is equal to 
$$
(1-u^2)^{\vert E \vert+ \vert V \vert}
\sum \varepsilon(\sigma) (-\alpha)^{\vert V_1\vert }
(-\beta)^{\vert V_2 \vert},
$$
extended over all sequences
$(V_1,V_2,V_3,\sigma,f,g)$. Now the last summation, say, $S$ is 
equal to
$$
\eqalignno{S&=\sum_{((V_1,V_2,V_3),\sigma)}
\varepsilon(\sigma)\,(-\alpha)^{|V_1|}(-\beta)^{|V_2|}
\prod_{i\in V_1} K(i,\sigma(i))\times \prod_{j\in V_2} \deg j\cr
&=\sum_{(V_1,V_2,V_3)}
\det(-\alpha\,K(i,j))_{(i,j\in V_1)}\times
(-\beta)^{|V_2|}\prod_{j\in V_2} \deg j \cr
&=\det(I-\beta(I+{\Cal Q})-\alpha\, {\Cal K}),\cr
}
$$
where $\Cal Q$ and $\Cal K$ are the two matrix ingredients of
the matrix $\Delta(u)$ defined in section~6.
Hence the sum of all the weights of the polite chaps
(and hence the sum of the weights of {\it all} chaps) equals
$$\displaylines{\qquad
(1-u^2)^{\vert E \vert+ \vert V \vert}
\det \Bigl(I+{u^2 \over 1-u^2}\,(I+{\Cal Q})-{{u} \over
{1-u^2}}\,{\Cal K}\Bigr)\hfill\cr
\kern4cm{} =
(1-u^2)^{\vert E \vert}
\det (I- u\,{\Cal K}+u^2\,{\Cal Q})\hfill\cr
\kern4cm{}=
(1-u^2)^{\vert E \vert}
\det \Delta (u),\hfill\cr}
$$
the right side of (7.1).

\subhead{\rm 7.3.}\enspace
Good and Bad Chaps\endsubhead
A chap is {\it hopelessly bad} if it contains superedges colored
$3a$, $3b$, $4a$, or $4b$. It is immediate that the sum of all the
weights of the hopelessly bad chaps is zero since superedges
colored $3a$ and $3b$ annihilate each other, as do those colored
$4a$ and $4b$. It is also clear that if a superedge~2 is present,
then the chap must be a hopelessly bad chap, since whenever a
vertex goes to an edge, some edge must go to a vertex through a
superedge necessarily colored $3a$, $3b$, $4a$, $4b$. For the
remaining chaps, the only way a vertex can be mapped onto is onto
itself (superedge~1), that explains the factor of
$(1-u^2)^{\vert V\vert}$ on the left side of (7.2). We can now
forget about the vertices and focus on the interaction of
the edges.

Having purged the hopelessly bad chaps, we can only have chaps
with superedges colored 5, 6, 7, that we shall further split into:

$e \buildrel 5a\over \longrightarrow e$ with weight 1;

$e \buildrel 5b\over \longrightarrow e$ with weight $-u^2$;

$e \buildrel 6a\over \longrightarrow e'$ 
($e'$ having the same origin as~$e$) with weight $u^2$;

$e \buildrel 6b\over \longrightarrow e'$ ($e'$ having the same
origin as~$e$ but $e'\not=e$) with weight
$u^2$;

$e\buildrel 7a\over \longrightarrow \overline e$ $(\overline
e=J(e))$ with weight $-u$;

$e\buildrel 7b\over \longrightarrow e'$ ($e'$ a successor
of~$e$ but different from $\overline e$) with weight $-u$.

\noindent
A not hopelessly bad chap is nevertheless {\it very bad}
if it contains superedges colored $5b$ or $6a$. These two cases
annihilate each other so we can easily execute all the very bad
chaps. It follows that a chap is {\it not very bad} if it contains
superedges colored $5a$, $6b$, $7a$, $7b$.

Finally, a chap Ch is a {\it bad
chap} if one of the three properties takes place:

(i) there is an edge $e$ such that 
$e\buildrel 5a\over \longrightarrow e$
and $\overline e\buildrel 6b\over \longrightarrow e'$
occur for some~$e'$;

(ii) there is an edge $e$ such that 
$e\buildrel 6b\over \longrightarrow e'$ and
$\overline e\buildrel 7b\over \longrightarrow e''$
occur for some~$e'$, $e''$ ;

(iii) there is an edge $e$ such that the sequence 
$e\buildrel 7a\over \longrightarrow \overline e
\buildrel 7b\over \longrightarrow e''$ occurs for some~$e''$.

\smallskip
If Ch is a bad chap let $e$ be the smallest offending edge.
We define ${\text {Ch}}'$ by making the obvious transposition, i.e., by
replacing the occurrence in case~(i) by the occurrence in case
(iii) and conversely, and replacing (ii) by
$e\buildrel 6b\over \longrightarrow e''$ and
$\overline e\buildrel 7b\over \longrightarrow  e'$.
It is clear that ${\text {Ch}}\mapsto {\text {Ch}}'$ is an involution of the set
of the not very bad chaps that reverses the sign and
preserves the absolute value of the weight. 

A non-bad chap will be called a {\it good chap}. It is then a chap
containing superedges colored $5a$, $6b$, $7a$, $7b$ and having
the following properties:

(i) whenever $e\buildrel 5a\over \longrightarrow e$
occurs, then either 
$\overline e\buildrel 5a\over \longrightarrow \overline e$, or
$\overline e\buildrel 7b\over \longrightarrow e'$ occurs;

(ii) whenever $e\buildrel 6b\over \longrightarrow e'$ occurs,
then either $\overline e\buildrel 7a\over \longrightarrow e$
or
$\overline e\buildrel 6b\over \longrightarrow e''$ occurs;

(iii) whenever  $e\buildrel 7a\over \longrightarrow \overline e$
occurs, then either
$\overline e\buildrel 6b\over \longrightarrow e'$ or
$\overline e\buildrel 7a\over \longrightarrow e$ occurs.

\subhead{\rm 7.4.}\enspace Enumerating the good chaps\endsubhead
Referring to the left side of (7.1) we are left to prove that the
weighted sum of all the good chaps is equal~to 
$$
\beta_u({\bold H})\times (1-u^2)^{|E|/2}
=\sum_{\pi\in {\Cal H}} (-1)^{\deg \pi}
u^{|\Cont(\pi)|}\times (1-u^2)^{|E|/2}.
$$
Say that $\tau$ is a {\it back-track involution}, if there exists
a subset $F(\tau)$ of the edge set~$E$ such that
$J(F(\tau))=F(\tau)$ and $\tau$ is the restriction of~$J$ to
$F(\tau)$. Let $\Cal T$ denote the set of all back-track
involutions. (Notice that $\Cal T\subset\Cal G$.) With those
notations we have
$$
(1-u^2)^{|E|/2}
=\sum_{\tau\in {\Cal T}} (-1)^{\deg \tau}u^{|F(\tau)|}.
$$
Denote by $w({\text {Ch}})$ and $\varepsilon({\text {Ch}})$ the weight
and the sign of a (good) chap~Ch, respectively. We are left to
prove the identity
$$
\sum_{{\text {Ch}}\ \text {good\ chap}} \kern-7pt 
\varepsilon({\text {Ch}})\,w({\text {Ch}})=
\sum_{\pi\in {\Cal H}} (-1)^{\deg \pi} u^{|\Cont(\pi)|}
\times
\sum_{\tau\in {\Cal T}} (-1)^{\deg \tau}u^{|F(\tau)|}.
$$

{\sl Construction of a bijection ${\text {Ch}}\mapsto (\pi,\tau)$ of
the set of good chaps onto ${\Cal H}\times {\Cal T}$}.\quad
The definition of a good chap shows that there are six cases to
consider depending on the colors of superedges going out of each
pair $e,\overline e$. The bijection is shown in the next table. 

\goodbreak
$$\hbox{\vbox{\offinterlineskip
\halign{\vrule\vrule depth 3pt height 12pt width
0pt\ \hfil$#$\hfil\  &\vrule\ \hfil$#$\hfil\ 
&\vrule\ \hfil$#$\hfil\ &\vrule\ \hfil$#$\hfil\ \vrule\cr
\noalign{\hrule}
&{\text {Ch}}&\pi&\tau\cr 
\noalign{\hrule}
(1)&e\buildrel 5a\over \longrightarrow e&&\cr
&\overline e\buildrel 5a\over \longrightarrow \overline e&&\cr
\noalign{\hrule}
(2)&e\buildrel 5a\over \longrightarrow e&\overline e
\longrightarrow e'&\cr
&\overline e\buildrel 7b\over \longrightarrow e'&&\cr
\noalign{\hrule}
(3)&e\buildrel 6b\over \longrightarrow e'&\overline e
\longrightarrow e'&e\longrightarrow \overline e\cr
&\overline e\buildrel 7a\over \longrightarrow e&&
\overline e\longrightarrow e\cr
\noalign{\hrule}}}\hskip.5cm
\vbox{\offinterlineskip\halign{\vrule\vrule depth 3pt height 12pt
width 0pt\ \hfil$#$\hfil\  &\vrule\ \hfil$#$\hfil\ 
&\vrule\ \hfil$#$\hfil\ 
&\vrule\ \hfil$#$\hfil\ \vrule\cr
\noalign{\hrule}
&{\text {Ch}}&\pi&\tau\cr 
\noalign{\hrule}
(4)&e\buildrel 6b\over \longrightarrow e'&e\longrightarrow e''&
e\longrightarrow \overline e\cr
&\overline e\buildrel 6b\over \longrightarrow e''&
\overline e\longrightarrow e'&\overline e\longrightarrow e\cr
\noalign{\hrule}
(5)&e\buildrel 7a\over \longrightarrow \overline e&& 
e\longrightarrow \overline e\cr
&\overline e\buildrel 7a\over \longrightarrow e&&
\overline e\longrightarrow e\cr
\noalign{\hrule}
(6)&e\buildrel 7b\over \longrightarrow e'&e
\longrightarrow e'&\cr
&\overline e\buildrel 7b\over \longrightarrow e''&
\overline e\longrightarrow e''&\cr
\noalign{\hrule}}}}
$$
For instance, in case (2) we define $\pi(\overline e)=e'$;
furthermore $e\not\in \Cont(\pi)$ and $e,\overline e\not\in
F(\tau)$. The definition of~$\tau$ is straightforward. 
To obtain~$\pi$ we start with the cycles of~Ch and make the
local modifications indicated. 

In the construction of~$\pi$ no edge $e$ is mapped onto its
reverse~$\overline e$, so that $\pi\in {\Cal H}$. The inverse
bijection is described by means of the same table.

What remains to be proved is the identity
$$\varepsilon({\text {Ch}})w({\text {Ch}})
=(-1)^{\deg\pi}\,u^{|\Cont(\pi)|}\,
(-1)^{\deg\tau}\,u^{|F(\tau)|}.\leqno(7.2)
$$
In cases~(1), (2) and~(6)
there is no modification in the composition of the cycles when
we go from~Ch to~$\pi$. In case~(3) the supervertex~$e$ is
deleted from the cycle containing~$\overline e$, but the
transposition $e\leftrightarrow \overline e$ occurs in~$\tau$.
In case~(4) two cycles of~$\pi$ are made out of a single one, or a
single cycle is made out of two existing ones. Therefore the sign
changes, but again $e\leftrightarrow \overline e$ occurs 
in~$\tau$. Finally, in case~(5) the transposition 
$e\buildrel 7b\over \longleftrightarrow \overline e$ is
transformed into the transposition $e\leftrightarrow \overline e$ 
in~$\tau$. Hence
$$\varepsilon(h)
= \varepsilon(\pi)\,\varepsilon(\tau)
$$
For each $i=1,\ldots,6$ let $n_i$ be the number of pairs
$(e,\overline e)$ falling into case~(i). 
The weight of~Ch (not counting the contribution due to
the vertices) is equal to
$$\eqalignno{w({\text {Ch}})
&=(-u)^{n_2}u^{2n_3}(-u)^{n_3}u^{4n_4}
(-u)^{2n_5}(-u)^{2n_6}\cr
&=(-u)^{n_2+n_3+2n_4+2n_6}\,
(-u)^{2n_3+2n_4+2n_5}\hfill\cr
&=(-u)^{|\Cont(\pi)|}\,(-u)^{|F(\tau)|}.\cr}
$$

\noindent
Altogether
$$
\eqalignno{
\varepsilon({\text {Ch}})w({\text {Ch}})
&=\varepsilon(\pi)\,(-u)^{|\Cont(\pi)|}\,
\varepsilon(\tau)\,(-u)^{|F(\tau)|}\cr
&=(-1)^{\deg\pi}\,u^{|\Cont(\pi)|}\,
(-1)^{\deg\tau}\,u^{|F(\tau)|}.\cr
}
$$

\head
8. A matrix-algebraic proof of (1.3)
\endhead
Let $(u(i,j))$ $(1\le i,j\le c_0)$ and $(v(i))$ $(1\le i\le c_0)$
be two sets of commuting variables. Introduce the {\it common
origin map} ``Com" as follows: if $e$ is an oriented edge that
goes from vertex~$i$ to vertex~$j$, define:
$$\leqalignno{
{\text {Com}}(e)&:=\sum_{e'\in E_i,\, e'\not=e} e';&(8.1)\cr
{\text {Com}}({\bold v})(e)&:=v(i)\,{\text {Com}}(e).&(8.2)\cr
}
$$
Thus ${\text {Com}}(e)$ is the sum of all edges, other than~$e$,
that have the same origin as~$e$. Keeping the same notations
we further define
$$\displaylines{\rlap{(8.3)}\hfill
{\text {Succ}}({\bold u})(e):=u(i,j)\,{\text {Succ}}(e),\hfill\cr
\noalign{\hbox{so that}}
\rlap{(8.4)}\hfill
A:=I+{\text {Succ}}({\bold u})+{\text {Com}}({\bold v})\hfill\cr}
$$
is an endomorphism of ${\Cal L}(E)$. Finally, for each $i=1,2,\ldots
,c_0$ let
$\Delta(i,i):=1+K(i,i)u(i,i)+ Q(i)v(i)$ and form the matrix
$$
\Delta=\left(\matrize{\Delta(1,1)&K(1,2)u(1,2)&\ldots
&K(1,c_0)u(1,c_0)\cr
K(1,2)u(1,2)&\Delta(2,2)&\ldots&K(2,c_0)u(2,c_0)\cr
\vdots&\vdots&\ddots&\vdots\cr
K(c_0,1)u(c_0,1)&K(c_0,2)u(c_0,2)&\ldots
&\Delta(c_0,c_0)\cr}\right).
$$

\proclaim{Proposition 8.1} The determinant of $A$ factorizes as
$$
\det A=\det \Delta \times\prod_{i=1}^{c_0} (1-v(i))^{Q(i)}.
\leqno(8.5)
$$
\endproclaim

\proof
There is no confusion in denoting both the
endomorphism and its corresponding matrix by the same symbol.
For each $i,j=1,2,\ldots,c_0$ let $A(i,j)$ be the linear map,
induced by $A$, that maps the space ${\Cal L}(E_j)$ into 
${\Cal L}(E_i)$. Its corresponding matrix is of dimention
$(Q(i)+1)\times (Q(j)+1)$. The matrix~$A$ itself is fully described by
the contents of all the blocks $A(i,j)$ $(i,j=1,2,\ldots,c_0)$. 

If $B$ is a matrix of order $n\times m$, denote by
$B_{1,\bullet}$, $B_{2,\bullet}$, \dots~, $B_{n,\bullet}$ its $n$ rows
(from top to bottom) and by
$B_{\bullet,1}$, $B_{\bullet,2}$, \dots~, $B_{\bullet,m}$ its $m$
columns (from left to right). Next define $\sigma \, B$ to be the matrix
whose rows are $B_{1,\bullet}$, $B_{2,\bullet}-B_{1,\bullet}$, \dots~,
$B_{n,\bullet}-B_{1,\bullet}$. Also let
$\alpha\,B$ be the matrix whose rows are
$B_{\bullet,1}+B_{\bullet,1}+\cdots+B_{\bullet,m}$,
$B_{\bullet,2}$, \dots~, $B_{\bullet,n}$.

First apply $\sigma\,\alpha$ to the
blocks $A(i,j)$ $(i>j)$ {\it below} the diagonal of the matrix~$A$
and  $\alpha\,\sigma$ to the other blocks $A(i,j)$ $(i\le j)$. It is
easily seen that those transformations keep invariant the value of
the determinant.  Its value does not change either if we further
make the following shift of rows {\it and} columns in the resulting
matrix:
$1\rightarrow 1$, $Q(1)+2\rightarrow 2$,
$Q(1)+Q(2)+3\rightarrow 3$, \dots~,
$Q(1)+\cdots+Q(c_0-1)+c_0\rightarrow c_0$.
We obtain the matrix
$$
D=\left(\matrize{\Delta&\star&\ldots&\star\cr
0&(1-v(1))I_{Q(1)}&\ldots&\star\cr
\vdots&\vdots&\ddots&\vdots\cr
0&0&\ldots&(1-v(c_0))I_{Q(c_0)}\cr}\right)
$$
where $\Delta$ is the matrix defined above.
Hence
$$
\det A=\det D=\det \Delta\times \prod_{i=1}^{c_0}
(1-v(i))^{Q(i)}.\qed 
$$

\medskip
Using the endomorphism $T={\text {Succ}}-J$ defined in section~6,
we have ${\text {Com}}=TJ$.
Accordingly, if we let $v(i)=u^2$ for all $i$ and replace all
the $u(i,j)$ by $-u$ in the definition of~$A$, we get
$A=I-u(T+J)+u^2TJ=(I-uT)(I-uJ)$. But
$\det(I-uJ )$ is clearly equal to $(1-u^2)^{c_1}$. Hence
$\det \Delta \prod\limits_{i=1}^n (1-u^2)^{Q(i)}=\det\Delta
(1-u^2)^{2c_1-c_0}=\det (1-uT)\det(I-uJ)$, so that
$\det(I-uT)=\det\Delta(1-u^2)^{c_1-c_0}$, which is Bass's
identity (1.3).

\bigskip\bigskip
\noindent
{\bf Acknowledgements}. The authors are thankful to Jean-Pierre
Jouanolou who convinced them that the first evaluation of the
Ihara-Selberg zeta function could be derived from Amitsur's identity
and to Jiang Zeng who drew their attentions to the paper by
Reutenauer and Sch\"utzenberger. They finally thank Phil Hanlon who
made several suggestions that improved the presentation of the paper.

\refstyle{C}

\enddocument


\Refs

\ref\no 1
\by Guido Ahumada
\paper Fonctions p\'eriodiques et formule des traces de Selberg sur
les arbres
\jour C. R. Acad. Sci. Paris
\vol 305
\yr 1987
\pages 709-712
\endref


\ref\no 2
\by S. A. Amitsur
\paper On the Characteristic Polynomial of a Sum of Matrices
\jour Linear and Multilinear Algebra
\vol 9
\yr 1980
\pages 177--182
\endref

\ref\no 3
\by Hyman Bass
\paper The Ihara-Selberg Zeta Function of a Tree Lattice
\jour Internat. J. Math.
\vol 3
\yr 1992
\pages 717-797
\endref


\ref\no 4
\by Pierre Cartier and Dominique Foata
\book Probl\`emes combinatoires de commutation et r\'ear\-ran\-ge\-ments
\publ Springer-Verlag (Lecture Notes in Math., vol.~85)
\publaddr Berlin
\yr 1969
\endref

\ref\no 5
\by K.T. Chen, R.H. Fox and R.C. Lyndon
\paper Free differential calculus, IV. The quotient groups of the lower
central series
\jour Ann. Math.
\vol 68
\yr 1958
\pages 81-95
\endref

\ref\no 6
\by Dominique Foata
\paper A combinatorial proof of Jacobi's identity
\jour Ann. Discrete Math.
\vol 6
\yr 1980
\pages 125-135
\endref

\ref\no 7
\by Jean-Pierre Jouanolou
\paper Personal communication
\yr 1996
\endref

\ref\no 8
\by M. Lothaire
\book Combinatorics on Words
\publ Addison-Wesley (Encyclopedia of Math. and its Appl., vol.~17)
\publaddr Reading, Mass.
\yr 1983
\endref

\ref\no 9
\by (Major) P.A. MacMahon
\book Combinatory Analysis, {\rm vol.~1}
\publ Cambridge Univ. Press (Reprinted by Chelsea, New York, 1955)
\publaddr Cambridge
\yr 1915
\endref

\ref\no 10
\by Sam Northshield
\paper  Proofs of Ihara's Theorem for Regular and
Irregular Graphs
\jour Proc. I.M.A. Workshop ``Emerging Applications
of Number Theory" (submitted)
\yr 1996
\endref

\ref\no 11
\by Dominique Perrin
\paper Personal communication
\yr 1996
\endref

\ref\no 12
\by Christophe Reutenauer and Marcel-Paul
Sch\"utzenberger
\paper A Formula for the Determinant of a Sum of
Matrices
\jour Letters of Math. Physics
\vol 13
\yr 1987
\pages 299--302
\endref

\ref\no 13
\by Gian-Carlo Rota
\paper On the Foundations of Combinatorial
Theory. I Theory of M\"obius Function
\jour Z.
Wahrscheinlichkeitstheorie
\vol 2
\yr 1964
\pages 340--368
\endref

\ref\no 14
\by Gian-Carlo Rota
\paper Report on the present state of
combinatorics {\rm (Inaugural address delivered at the 5th Formal
Power Series and Algebraic Combinatorics Conference,
Florence, 21 June 1993)}
\jour Discrete Math.
\vol 153
\yr 1996
\pages 289--303
\endref

\ref\no 15
\by Marcel-Paul Sch\"utzenberger
\paper Sur une propri\'et\'e
combinatoire des alg\`ebres de Lie libres pouvant \^etre utilis\'ee
dans un probl\`eme de math\'ematiques appliqu\'ees
\inbook S\'eminaire d'alg\`ebre et de th\'eorie des nombres 
[P. Dubreil, M.-L. Dubreil-Jacotin, C. Pisot, 
1958-59]
\publ Secr\'etariat Math\'ematique
\publaddr 11, rue Pierre-Curie, F-75005
\yr 1960
\pages 1-01--1-13
\endref

\ref\no 16
\by Marcel-Paul Sch\"utzenberger
\paper On a factorization
of free monoids
\jour Proc. Amer. Math. Soc.
\vol 16
\yr 1965
\pages 21-24
\endref

\ref\no 17
\by Richard P. Stanley
\book Enumerative Combinatorics, {\rm vol.
1}
\publ Wadsworth \& Brooks
\publaddr Monterey
\yr 1986
\endref

\ref\no 18
\by H. M. Stark and A. A. Terras
\paper Zeta Functions of Finite Graphs and Coverings
\jour Adv. in Math.
\vol 121
\yr 1996
\pages 124-165
\endref


\ref\no 19
\by G\'erard Viennot
\book Alg\`ebres de Lie libres et mono\"\i des libres
\publ Springer-Verlag (Lecture Notes in Math., vol.~691)
\publaddr Berlin
\yr 1978
\endref

\ref\no 20
\by Doron Zeilberger
\paper A Combinatorial Approach to
Matrix Algebra
\jour Discrete Math.
\vol 56
\yr 1985
\pages 61-72
\endref

\endRefs
\bye